\documentclass{article}
\usepackage{amssymb,amsfonts,amsmath,theorem,amscd}

\newcommand{\FF}{\mathbb F}

\newcommand{\ZZ}{\mathbb Z}

\newtheorem{theorem}{Theorem}
\newtheorem{lemma}[theorem]{Lemma}

{\theorembodyfont{\rmfamily}}

\begin{document}

\title{The distribution of the number of points modulo an integer on elliptic curves over finite fields}
%\titlerunning{The distribution of $\#E(\mathbb{F}_q)$ modulo $N$}
\author{Wouter Castryck\thanks{
\emph{Address:} Katholieke Universiteit Leuven,
Departement Wiskunde,
Celestijnenlaan 200B, 3001 Leuven (Heverlee), Belgium;
\emph{E-mail:} firstname.lastname@wis.kuleuven.be; 
%\emph{Tel:} (+32) 16 32 70 06; \emph{Fax:} (+32) 16 32 79 98
}\ and Hendrik Hubrechts${}^\ast$}
\date{}
\maketitle

\begin{abstract}
Let $\mathbb{F}_q$ be a finite field and let $b$ and $N$ be
integers. We prove explicit estimates for the probability
that the number of rational points on a
randomly chosen elliptic curve $E$ over $\mathbb{F}_q$ equals $b$
modulo $N$. The underlying tool is an equidistribution result
on the action of Frobenius on the $N$-torsion subgroup of $E$.
Our results subsume and extend previous work by Achter and Gekeler.\\

\noindent Keywords: elliptic curves, finite fields, Frobenius statistics, modular curves\\

\noindent MSC 2010: 14H52, 14K10\\
\end{abstract}

\section{Introduction}

If one writes the number of rational points on an elliptic curve
$E$ over a finite field $\mathbb{F}_q$ as $q + 1 - T_E$, then the
integer $T_E$ is called the trace of Frobenius of $E$. Hasse proved
that $T_E \in [-2\sqrt{q}, 2\sqrt{q}]$, but within this interval the
trace of Frobenius is an unpredictable number, seemingly picked at
random. Since the 1960's, its statistical behavior has become
subject to extensive study.

To make the problem well-defined, the best-known approach is to
fix an elliptic curve $E$ over a number field $K$ and to consider
it modulo various prime ideals $\mathbf{p} \subset \mathcal{O}_K$
of good reduction. Based on experimental evidence, Sato and Tate
conjecturally described how the traces of Frobenius of $E$ mod
$\mathbf{p}$ are (after being normalized by
$2\sqrt{N(\mathbf{p})}$) distributed along the interval $[-1,1]$. We refer to \cite{Carayol} for the details
and an introduction to the recent progress on this matter.
%To each
%$\mathbf{p}$ one can associate a number $t_\mathbf{p} \in [-1,1]$, obtained
%by dividing the trace of Frobenius of $E/\mathbf{p}$ by $2\sqrt{N(\mathbf{p})}$.
%Sato and Tate independently conjectured that, when $E$ has no complex multiplication,
%the numbers $t_\mathbf{p}$ are randomly distributed according to the semicircular
%probability measure $\mu = \frac{2}{\pi}\sqrt{1 - t^2} dt$.
%More precisely, for each measurable set $S \subset [-1,1]$ the limit
%\[ \lim_{n \rightarrow \infty} \frac{\# \left\{ \left. \mathbf{p} \subset \mathcal{O}_K \, \right|
%\, N(\mathbf{p}) \leq n \text{ and } t_\mathbf{p} \in S \right\} }{\# \left\{ \left. \mathbf{p} \subset \mathcal{O}_K \, \right|
%\, N(\mathbf{p}) \leq n\right\}}\]
%exists and is equal to $\mu(S)$.
%Though the conjecture remains unproven, recently it has been
%virtually settled for totally real number fields by (mainly)
%Taylor \cite{Carayol}.
%If $E$ does have complex multiplication, work by Hecke implies that the $t_\mathbf{p}$ are randomly
%distributed with respect to $dt/(\pi \sqrt{1 - t^2})$.
%See \cite{Carayol} for a detailed introduction
%to the Sato-Tate conjecture.

Another approach is to fix the finite field $\mathbb{F}_q$ and to consider all $\mathbb{F}_q$-isomor\-phism
classes of elliptic curves $E$ over it. Their traces of Frobenius $T_E$
define a discrete probability measure $\mu_q$ on $\{ -\lfloor2\sqrt{q}\rfloor, \dots, \lfloor2\sqrt{q}\rfloor \}$.
%Since over number fields, curves with complex multiplication are comparatively rare,
%a natural guess is that $\mu_q \rightarrow \mu$ as $q$ gets big.
As above, one can normalize to obtain a distribution $\tilde{\mu}_q$ on $[-1,1]$.
Birch \cite{Birch}, Deligne \cite[3.5.7]{Deligne} and Yoshida \cite{Yoshida} proved results on the limiting
behavior of $\tilde{\mu}_q$ as $q$ tends to infinity, thereby
lending indirect support for the Sato-Tate conjecture.
%One can prove that $\mu_q \rightarrow \mu$ as $q$ tends to
%infinity. This result has been obtained in part by Deligne who
%studied the limit behavior of $\mu_{q^k}$ with $q$ fixed and $k$
%tending to infinity \cite[3.5.7]{Deligne}, and in part by Birch who dealt
%with $\mu_{p^k}$, with $k$ fixed and the prime $p$ tending to
%infinity \cite{Birch}. VERIFI\"EREN.
However, some interesting properties that are related to the discrete
nature of $\mu_q$ become dissolved in the limit procedure.
%since the discrete measures $\mu_q$ have
%many intriguing properties that become dissolved in taking the
%limit $\mu_q \rightarrow \mu$.
As an introductory exercise, the reader
is invited to show that when $q$ is odd, $T_E$ favors even numbers. This
is related to the fact that a randomly chosen cubic polynomial
$f(x) \in \mathbb{F}_q[x]$ has a rational root with a probability
that tends to $\frac{2}{3}$ as $q$ tends to infinity. More generally, for any integer
$N \geq 2$, the probability that $\#E(\mathbb{F}_q) = q + 1 - T_E$ is
divisible by $N$ tends to be strictly bigger than $\frac{1}{N}$.
Lenstra was the first to observe this phenomenon in his famous paper
\cite[Prop.~1.14]{Lenstra} --- which
has implications for integer factorization \cite{Lenstra} and
cryptography \cite{GalbraithMcKee} --- and he proved explicit
estimates in the situation where $N$ is a prime number different from $p = \text{char}(\FF_q)$.
The proof
uses modular curves and
was generalized to arbitrary $N$ by Howe \cite[Thm.~1.1]{Howe}.

In this article, we study the more general question of how $\#E(\mathbb{F}_q)$ modulo $N$ is distributed
along $\{0, 1, \dots, N-1\}$.
%For any positive integer $N$, the traces $T_E$ mod $N$ define a probability measure on $\{0, \dots, N-1 \}$.
For an arbitrary integer $N\geq 2$ and $t \in \mathbb{Z}_{\geq 0}$, write $P_{q,N}(t)$ for the probability that $T_E \equiv t \bmod N$.
We prove:

\begin{theorem} \label{maintheorem}
Let $Q = \{ \, p^k \, | \, p \text{ prime}, k \in \mathbb{Z}_{\geq 1} \, \}$, and let
$r : Q \times \mathbb{Z}_{\geq 1} \times \mathbb{Z}_{\geq 0} \rightarrow \mathbb{Q}_{\geq 0}$ be the unique function
satisfying the following rules.
\begin{itemize}
\item[(i)] $r$ is a multiplicative arithmetic function in the second argument, i.e.\ for all $q \in Q$,
$t \in \mathbb{Z}_{\geq 0}$ and coprime $M,N \in \mathbb{Z}_{\geq 1}$ one has
\[ r(q,MN,t) = r(q,M,t) \cdot r(q,N,t). \]
\item[(ii)] If $N = \ell^n$ for an integer $n \geq 1$ and a prime number $\ell$, then for all $q \in Q$ that are coprime to $\ell$
     and all $t \in \mathbb{Z}_{\geq 0}$ one has
\[ r(q,N,t) = \frac{\Psi(t^2 - 4q)}{\ell^{3n} - \ell^{3n-2}} \]
for the function $\Psi : \mathbb{Z} \rightarrow \mathbb{Z}$ that
is described explicitly in Section~\ref{traceNdist} below.
In case $\ell \geq 3$ and $n=1$ we have $\Psi : x \mapsto \ell^2 + \left( \frac{x}{\ell} \right) \ell$, where
$\left( \frac{\cdot}{\cdot} \right)$ is the Legendre symbol.
\item[(iii)] If $N=p^e$ for an integer $e \geq 1$ and a prime number $p$, then for all $q \in Q$ that are a power of $p$
    and all $t \in \mathbb{Z}_{\geq 0}$ one has
\[ r(q,N,t) = \begin{cases}
\displaystyle \frac{1}{p^e - p^{e-1}} & \textnormal{if \ } t \not \equiv 0 \bmod p, \\
\displaystyle 0 & \textnormal{if \ } t \equiv 0 \bmod p.  \end{cases} \]
\end{itemize}
Then there exists an absolute and explicitly computable constant $C \in \mathbb{R}_{>0}$ such that
for all $q \in Q, N \in \mathbb{Z}_{\geq 1}$ and $t \in \mathbb{Z}_{\geq 0}$ one has
\[ \left| P_{q,N}(t) - r(q,N,t) \right| \leq C \cdot \frac{N^2 \ln \ln N}{\sqrt{q}}.\]
\end{theorem}
Theorem~\ref{maintheorem} is essentially obtained from an equidistribution
result on the action of
$q$th power Frobenius on the $N$-torsion group
$E[N]$ of $E$. Throughout this article, for any integer $A$ we will
write $Z_A$ for the ring of residues $\mathbb{Z}/(A)$. Factor $N$ as $N'p^e$, where $N'$ is coprime to $p$.
In case $e \geq 1$, we suppose
that $E$ is taken from the set of
ordinary elliptic curves. Then $E[N] \cong E[N'] \oplus E[p^e] \cong Z_{N'} \oplus Z_{N'}
\oplus Z_{p^e}$.
With respect to a $Z_{N'}$-module basis of $E[N']$ and a generator of $E[p^e]$,
the action of $q$th power Frobenius is given by a pair
\[ (F,T) \in \text{GL}_2(Z_{N'}) \oplus Z_{p^e}^\times\]
satisfying $\det F \equiv q \bmod N'$, $\text{Tr} \, F \equiv T_E \bmod N'$, and $T \equiv T_E \bmod p^e$.
When considering all bases of $E[N']$, the corresponding matrices $F$ yield a conjugacy class of $\text{GL}_2(Z_{N'})$ which we denote
by $\mathcal{F}_E$. In contrast, the element $T$ does not depend
on the generator of $E[p^e]$; for sake of consistency, we will denote it by $\mathcal{T}_E$.
Note that $(\mathcal{F}_E,\mathcal{T}_E)$ can in fact be seen as a conjugacy class of $\text{GL}_2(Z_{N'}) \oplus Z_{p^e}^\times$.
Denote the subset of $\text{GL}_2(Z_{N'})$
consisting of all matrices of determinant $q$ by $\mathcal{M}_{q,N'}$.
Then the equidistribution theorem reads:
\begin{theorem} \label{equidistributionmain}
There exists an absolute and explicitly computable constant $C \in \mathbb{R}_{>0}$ such that for
every conjugacy class $\mathcal{F} \subset \emph{GL}_2(Z_{N'})$ of
matrices of determinant $q$ and every element $\mathcal{T} \in Z_{p^e}^\times$ one has
\[ \left| P_\mathcal{F,\mathcal{T}} - \frac{\#\mathcal{F}}{\#\mathcal{M}_{q,N'}} \cdot \frac{1}{\varphi(p^e)} \right| \leq C \cdot
\frac{p^eN'^2 \ln \ln N'}{\sqrt{q}},\] where $\varphi$ is Euler's totient function and
$P_{\mathcal{F},\mathcal{T}}$ is the probability that
\begin{itemize}
\item[-] $\mathcal{F}_E = \mathcal{F}$ if $e=0$;
\item[-] $E$ is ordinary and $(\mathcal{F}_E, \mathcal{T}_E) = (\mathcal{F},\mathcal{T})$ if $e \geq 1$.
\end{itemize}
\end{theorem}
Loosely stated: if $q$ gets big, a Frobenius conjugacy class
becomes as likely as its own relative size.

Theorem~\ref{equidistributionmain} fits within the random matrix philosophy that dominates nowadays
research on the statistical behavior of Frobenius, both in the
Sato-Tate setting (fixed curve, varying field) as in the
setting of a fixed field and a varying curve. This was
initialized by Deligne, who obtained his earlier-mentioned result as
a consequence to an equidistribution theorem in \'etale
cohomology. The random matrix idea has proven to provide
well-working models for higher genus analogues of the Frobenius
distribution problem \cite{KatzSarnak,KedlayaSutherland}, although
many statements remain conjectural. The standard reference has
become the book by Katz and Sarnak \cite{KatzSarnak}, who also
refined Deligne's equidistribution theorem to a version
\cite[9.7]{KatzSarnak} that was used by Achter to prove a variant
of Theorem~\ref{equidistributionmain} that works in arbitrary
genus \cite[Thm.~3.1]{Achter08}. However, Achter's result involves
certain mild restrictions on $q$
and $N$, the main one
being that $q$ and $N$ should be coprime. Our attention will be devoted to a slightly more elementary
approach, based on the modular covering $X(p^e; \zeta_{N'}) \rightarrow X(1; 1)$
and Chebotarev's density theorem for
function fields. Apart from resolving the conditions on $q$ and $N$,
this has the additional advantage of yielding a tighter error bound: in Achter's case
it is of the form $C \cdot N^3 / \sqrt{q}$. There is no doubt that
several specialists in the field expected an approach using Chebotarev's density theorem to work, but up to our knowledge, a complete
proof of Theorem~\ref{equidistributionmain}
did not appear in the literature before.

Given Theorem~\ref{equidistributionmain}, the proof of Theorem~\ref{maintheorem}
then comes down to determining the number of matrices in $\mathcal{M}_{q,N'}$ with trace $t$.
This is elaborated in Section~\ref{traceNdist}.
Again, large parts of this matrix count
have been carried out before, now by Gekeler \cite{Gekeler03}, who
worked towards estimating $P_{q,N}(t)$ under
certain mild conditions on $q$ and $N$, while assuming the random matrix principle as a black box.
The current article can therefore be viewed as a self-contained subsumption of both Achter's result in genus $1$ and Gekeler's count,
providing more elementarily flavored proofs and removing the restrictions on $q$ and $N$.

It is worth noting that Theorem~\ref{equidistributionmain} can be used
to study a number of alternative questions, various of which have been addressed in the literature before, albeit often
conjecturally. E.g.,
in the weaker set-up where $\mathbb{F}_q$
is a large prime field that is chosen at random, Gekeler
studied the probability that $E[\ell^\infty](\mathbb{F}_q)$ has a
given structure and the probability that $E(\mathbb{F}_q)$ is cyclic
\cite{Gekeler06,Gekeler08}. The latter probability has also been
investigated by Vl\v{a}du\c{t} in case $\mathbb{F}_q$ is fixed
\cite{Vladut}, building on Howe's aforementioned work. Still for $\mathbb{F}_q$ fixed,
Galbraith and McKee conjecturally estimated the chance that
$E(\mathbb{F}_q)$ is a prime number \cite{GalbraithMcKee}.
Achter
and Sadornil studied the probability that $E$ has a given number of
rational isogenies of given prime degree emanating from it
\cite{AchterSadornil}. For higher genus curves $C / \mathbb{F}_q$,
Achter gave explicit estimates for the chance that
$\text{Jac}(C)[N](\mathbb{F}_q)$ has a given structure
\cite{Achter06,Achter08}, and Chavdarov proved that the numerator
of the zeta function $Z_C(T)$ is generically irreducible
\cite{Chavdarov}. Recently, the current authors, Folsom and Sutherland \cite{genus2prime}
studied the probabilities of having prime order
and of cyclicity of $\text{Jac}(C)[N](\mathbb{F}_q)$.

One interesting question that did not see explicit study so far is
on the probability $P'_q(N)$ that $E$ contains a rational point of given order $N$.
In Section~\ref{sectionOrder} we prove:
\begin{theorem} \label{mainthmorder}
Let $r' : Q \times \mathbb{Z}_{\geq 1} \rightarrow \mathbb{Q}_{\geq 0}$ be the unique
function satisfying the following rules.
\begin{itemize}
  \item[(i)] $r'$ is a multiplicative arithmetic function in the second argument, i.e.\ for
  all $q \in Q$ and coprime $M,N \in \mathbb{Z}_{\geq 1}$ one has
  \[ r'(q,MN) = r'(q,M) \cdot r'(q,N). \]
  \item[(ii)] If $N = \ell^n$ for an integer $n \geq 1$ and a prime number $\ell$, then for all  $q \in Q$ that are coprime to $\ell$
      one has
  \[ r'(q,N) = \left\{ \begin{array}{ll} 1/(\ell^n - \ell^{n-2}) & \text{if $\nu \geq n$}, \\ (\ell^{2 \nu + 1} + 1)/(\ell^{n + 2\nu - 1}(\ell^2-1)) & \text{if $\nu < n$}, \\ \end{array} \right. \]
  where $\nu$ is the $\ell$-adic valuation of $(q-1)$.
  \item[(iii)] If $N = p^e$ for an integer $e \geq 1$ and a prime number $p$, then for all $q \in Q$ that are a power of $p$
      one has
  \[ r'(q,N) = 1/(p^e - p^{e-1}).\]
\end{itemize}
Then there exists an absolute and explicitly computable constant $C \in \mathbb{R}_{> 0}$ such that
for all $q \in Q$ and $N \in \mathbb{Z}_{\geq 1}$ one has
\[ \left| P'_q(N) - r'(q,N) \right| \leq C \cdot \frac{N^2 \ln \ln N}{\sqrt{q}}.\]
\end{theorem}
In fact,
this theorem can be derived directly from the work of Howe \cite[Thm.~1.1]{Howe}, instead of
Theorem~\ref{equidistributionmain}.

\section{Equidistribution of Frobenius} \label{matrixdist}

%This section briefly discusses some features of the
%modular curve $X(N)$ and the way it covers $X(1) \cong
%\mathbb{P}^1$.
%An implicit reference for this section are the lecture notes by
%Deligne and Rapoport \cite{DeligneRapoport} and the earlier
%work by Igusa \cite{Igusa1,Igusa2} on which these build.

In this section, we will prove Theorem~\ref{equidistributionmain}.
The two main
theoretical ingredients are the modular curve $X(p^e; \zeta_{N'})$ and Chebotarev's density theorem
for function fields.

We first recall some facts on modular curves.
Let $\FF_q$ be a finite field of characteristic $p$ having $q$ elements, let
$N$ be a positive integer, and write $N = N'p^e$ with $N'$ coprime to $p$. Assume throughout that $N' > 2$. Fix
an algebraic closure $\overline{\mathbb{F}}_q$ of $\mathbb{F}_q$ and a primitive
$N'$th-root of unity $\zeta_{N'} \in \overline{\FF}_q$.
Consider all quartets $(E,P,Q,R)$ for which $E$ is an elliptic
curve over $\overline{\mathbb{F}}_q$ and
\begin{itemize}
  \item[-] $P, Q \in E[N']$ satisfy
$e_{N'}(P,Q) = \zeta_{N'}$, where
\[ e_{N'} : E[N'] \times E[N'] \rightarrow \{\text{$N'$th-roots of unity} \} \]
is the Weil pairing \cite[III.\S 8]{Silverman};
  \item[-] $R \in E^{(p^e)}$
is a generator of the kernel of the
Verschiebung $V_{p^e} : E^{(p^e)} \rightarrow E$,
where $E^{(p^e)}$ is the elliptic curve obtained by raising all coefficients of a model of $E$ to
the $p^e$th power.
\end{itemize}
Two
quartets $(E,P,Q,R)$ and $(E',P',Q',R')$ are called equivalent
if there exists an $\overline{\mathbb{F}}_q$-isomorphism $E
\rightarrow E'$ mapping $P$ to $P'$ and $Q$ to $Q'$, such that
the induced isomorphism $E^{(p^e)} \rightarrow E'^{(p^e)}$ takes
$R$ to $R'$. As a special
instance, using multiplication by $-1$, we have that $(E,P,Q,R)$ is
equivalent to $(E,-P,-Q,-R)$. Denote
the set of equivalence classes of such quartets by $Y(p^e; \zeta_{N'})$.
Then there exists an irreducible nonsingular projective curve
$X(p^e;\zeta_{N'})$ over $\overline{\mathbb{F}}_q$, along with a morphism
\[ J : X(p^e; \zeta_{N'}) \rightarrow \mathbb{P}^1 \supset \text{Spec} \, \overline{\mathbb{F}}_q[j] \]
such that
\begin{itemize}
\item[-] the points of $J^{-1}(\text{Spec} \, \overline{\mathbb{F}}_q[j])$
are in bijective correspondence with $Y(p^e; \zeta_{N'})$,
giving the latter the structure of an irreducible nonsingular affine curve over $\overline{\mathbb{F}}_q$;
\item[-] if $x \in J^{-1}(\text{Spec} \, \overline{\mathbb{F}}_q[j])$
corresponds to a quartet $(E,P,Q,R)$, then $J(x) = j(E)$, the $j$-invariant of $E$;
\item[-] $J$ is a Galois covering with Galois group $\left( \text{SL}_2(Z_{N'}) \oplus Z_{p^e}^\times \right)/ \{ \pm 1 \}$,
where $\{ \pm 1 \}$ is understood to be diagonally embedded;
the action is such that an element
\[ \pm \left( \begin{pmatrix} a & b \\ c & d \\ \end{pmatrix}, u \right) \]
takes the point corresponding to the class $(E,P,Q,R) \in Y(p^e; \zeta_{N'})$ to
the point corresponding to $(E,aP + cQ, bP + dQ, uR)$.
\end{itemize}
Moreover, $X(p^e; \zeta_{N'})$ and $J$ are naturally defined over
$\mathbb{F}_q(\zeta_{N'})$, and for $k = [\mathbb{F}_q(\zeta_{N'}) : \mathbb{F}_q]$
the action of $q^k$th power Frobenius on $X(p^e; \zeta_{N'})$ commutes with
the action of the Galois group.
When restricted to $Y(p^e; \zeta_{N'})$, the Frobenius action
is given by $(E,P,Q,R) \mapsto (E^{(q^k)},P^{(q^k)},Q^{(q^k)},R^{(q^k)})$.
The genus
of $X(p^e; \zeta_{N'})$ equals
\begin{equation} \label{modulargenus}
 \left\{ \begin{array}{ll} 1 + \frac{1}{24}(N-6) \varphi(N) \tilde{\varphi}(N) & \text{if $e = 0$} \\
 1 + \frac{1}{48}(N-12)\varphi(N) \tilde{\varphi}(N') & \text{if $e \geq 1$}. \\ \end{array} \right.
\end{equation}
%while the number of $x \in Y(p^e;\zeta_{N'})$ for which $J(x)$ is supersingular
%is given by
%\begin{equation} \label{supersingularlocus}
%\frac{p-1}{24} N' \varphi(N') \psi(N').
%\end{equation}
Here $\varphi : x \mapsto x \prod_{p \mid x}(1 - 1/p)$ is Euler's totient function,
and $\tilde{\varphi}$ is, somehow dually, defined by $x \mapsto x \prod_{p \mid x}(1 + 1/p)$.

For a proof of the above statements, we refer to the notes by Deligne and Rapoport \cite{DeligneRapoport}, the article of Howe \cite[Prop.~3.1~and~3.2]{Howe}
and the many references therein to the book of Katz and Mazur \cite{KatzMazur}.
In the latter, the curve $X(p^e; \zeta_{N'})$ is denoted $\overline{\mathfrak{M}}(\mathcal{P})$,
where $\mathcal{P}$ is the moduli problem $([\Gamma(N')]^\text{can},[\text{Ig}(p^e)])$ over $(\text{Ell}/ \overline{\mathbb{F}}_q)$.
Howe denotes
this curve by $\overline{X}(N',N)$.
The condition $N' > 2$ is necessary
for $[\Gamma(N')]^\text{can}$ to be representable in the sense of \cite[(4.3)]{KatzMazur}; see also
\cite[(10.9.3)]{KatzMazur}. It is possible to construct similar modular curves for $N' \leq 2$, as
illustrated by Howe \cite[Prop.~3.1]{Howe}, but we will not need this.

The $\mathbb{F}_q(\zeta_{N'})$-rational morphism $J$
gives rise
to a field extension
\[ \mathbb{F}_q(\zeta_{N'})(j) \ \subset \ \mathbb{F}_q(\zeta_{N'})(X(p^e; \zeta_{N'})) =: L. \]
Our central object of interest will be the larger extension
\[ K := \mathbb{F}_q(j) \ \subset \ \mathbb{F}_q(\zeta_{N'})(j) \ \subset \ L.\]
It allows a modular interpretation as follows.
Let $H \subset Z_{N'}^\times$ be the group generated by $q$ mod $N'$.
Under the map
$h \mapsto \zeta_{N'}^h$, its elements are
in bijective correspondence with the $\text{Gal}(\mathbb{F}_q(\zeta_{N'}), \mathbb{F}_q)$-orbit
of $\zeta_{N'}$.
Then similar to before, using $[\Gamma(N')]^{\mathbb{Z}[\zeta_{N'}]^H\text{-can}}$ instead of $[\Gamma(N')]^\text{can}$, we can define a complete nonsingular (but possibly reducible) curve $X^H(p^e; \zeta_{N'})$ over
$\overline{\mathbb{F}}_q$ along with a morphism $J^H$ to $\mathbb{P}^1$,
such that
\begin{itemize}
\item[-]
$(J^H)^{-1}
(\text{Spec} \, \overline{\mathbb{F}}_q[j])$ can be identified with
\[
  Y^H(p^e; \zeta_{N'}) \ := \ Y(p^e; \zeta_{N'}) \ \sqcup \ Y(p^e; \zeta_{N'}^q) \ \sqcup \ \dots \ \sqcup \ Y(p^e; \zeta_{N'}^{q^{k-1}});\]
\item[-] if $x \in (J^H)^{-1}(\text{Spec} \, \overline{\mathbb{F}}_q[j])$
corresponds to a quartet $(E,P,Q,R)$, then $J^H(x) = j(E)$;
\item[-] $J^H$ is a Galois covering with Galois group
$G = \left(\text{GL}_2^H(Z_{N'}) \oplus Z_{p^e}^\times \right) / \{ \pm 1 \}$,
where $\text{GL}_2^H(Z_{N'})$ is the group of matrices of $\text{GL}_2(Z_{N'})$ taking
determinants in $H$; on $Y^H(p^e; \zeta_{N'})$, the action of an element
\begin{equation} \label{galoisaction}
  \pm \left(\begin{pmatrix} a & b \\ c & d \\ \end{pmatrix}, u \right)
\end{equation}
is such that it
takes the point corresponding to the class $(E,P,Q,R)$ to
the point corresponding to $(E,aP + cQ, bP + dQ, uR)$.
\end{itemize}
Moreover, $X^H(p^e; \zeta_{N'})$ and $J^H$ can be defined over $\mathbb{F}_q$, and the action
of $q$th power Frobenius on $Y^H(p^e; \zeta_{N'})$ is given by
$(E,R,P,Q) \mapsto (E^{(q)},R^{(q)},P^{(q)},Q^{(q)})$. This action commutes
with the action of the Galois group. Considered as a scheme over $\mathbb{F}_q$, the curve $X^H(N'; \zeta_{N'})$
is irreducible, and the function field extension corresponding to the rational
morphism $J^H$ to $\mathbb{P}^1$ is nothing else than $K \subset L$. In particular, this extension is Galois.

We now start working towards an application of Chebotarev's density theorem
to $K \subset L$. Let $R = \mathbb{F}_q[j]$ and let $S$ be its integral closure inside $L$.
Then $S$ is the coordinate ring of
$Y^H(p^e; \zeta_{N'})$ considered as an $\mathbb{F}_q$-scheme.
Let $j_0 \in \FF_q$ and let $E_{j_0} / \mathbb{F}_q$
be an elliptic curve with $j$-invariant $j_0$. A quartet $\mathcal{E} = (E_{j_0},P,Q,R)$ on
$Y^H(p^e; \zeta_{N'})$ corresponds to a maximal ideal
$\mathbf{m}_\mathcal{E}$ in $S\otimes \overline{\FF}_q$. Define
$\mathbf{P}_\mathcal{E} := \mathbf{m}_\mathcal{E}\cap S$, which
can be viewed as a closed point of $Y^H(p^e; \zeta_{N'})$ as an $\mathbb{F}_q$-scheme.
Suppose that $\mathbf{P}_\mathcal{E}$ is unramified over $K$,
which is equivalent to
\begin{equation} \label{unramifiedcondition}
 j_0 \ \not \in \left\{ \begin{array}{ll} \{0, 1728\} & \text{if $e = 0$} \\
  \mathcal{J}_\text{ss} \cup \{0,1728\} & \text{if $e \geq 1$}, \\ \end{array} \right.
\end{equation}
where $\mathcal{J}_\text{ss} \subset \overline{\mathbb{F}}_q$ is the set
of supersingular $j$-invariants.
As explained in \cite[Section 6.2]{Fried} we can
associate to $\mathbf{P}_\mathcal{E}$ its Frobenius
automorphism
$\left[\frac{L/K}{\mathbf{P}_\mathcal{E}}\right]\in G = \text{Gal}(L/K)$.
With $\mathbf{p}_\mathcal{E}:=\mathbf{P}_\mathcal{E}\cap R$, this
automorphism is uniquely determined by the condition
\begin{equation}\label{eqFrobCongruence}\left[\frac{L/K}{\mathbf{P}_\mathcal{E}}\right] x \equiv x^{N(\mathbf{p}_\mathcal{E})}\bmod \mathbf{P}_\mathcal{E},\qquad\text{for all }x\in S.\end{equation}
We note that $j_0 \in\FF_q$ implies that $\mathbf{p}_\mathcal{E} =
(j-j_0)$ and hence $N(\mathbf{p}_\mathcal{E}) = q$.
Geometrically, condition \eqref{eqFrobCongruence} just means that if
\[ \{(E_{j_0},P_1,Q_1,R_1), (E_{j_0},P_2,Q_2,R_2), \dots, (E_{j_0},P_{\deg \mathbf{P}_\mathcal{E}},Q_{\deg \mathbf{P}_\mathcal{E}},R_{\deg \mathbf{P}_\mathcal{E}}) \}\]
is the set of points of $Y^H(p^e;\zeta_{N'})$ (maximal ideals of $S \otimes
\overline{\FF}_q$) above $\mathbf{P}_\mathcal{E}$, then
$\left[\frac{L/K}{\mathbf{P}_\mathcal{E}}\right] \in
G$ permutes this set in the same manner
as described in (\ref{galoisaction}) above. If $\mathbf{P}'$ is another prime ideal of $S$
above $\mathbf{p}_{\mathcal{E}}$, we have that the Frobenius
automorphism $\left[\frac{L/K}{\mathbf{P}'}\right]$ is conjugated
to $\left[\frac{L/K}{\mathbf{P}_{\mathcal{E}}}\right]$. The
Artin symbol
\[\left(\frac{L/K}{\mathbf{p}_\mathcal{E}}\right)\]
of $\mathbf{p}_\mathcal{E}$ is then defined as the conjugacy class
of $\left[\frac{L/K}{\mathbf{P}_{\mathcal{E}}}\right]$ in
$G$. Thus, the Artin symbol associated to a prime ideal
$\mathbf{p}_\mathcal{E} = (j-j_0) \subset R$ with $j_0$ satisfying (\ref{unramifiedcondition})
is the conjugacy class of $G = \left(\text{GL}_2^H(Z_{N'}) \oplus Z_{p^e}^\times \right) / \{ \pm 1 \}$
obtained by considering for an elliptic curve $E/ \mathbb{F}_q$ with $j$-invariant $j_0$
the action of $q$th power Frobenius with respect to
\begin{itemize}
\item[-] all generators $R$ of $\text{ker} \, V_{p^e}$;
\item[-] all bases $P,Q$ of $E[N']$ for which $e_{N'}(P,Q) = \zeta_{N'}^h$ for some $h \in H$.
\end{itemize}

Chebotarev's density theorem, in the following version of Fried and Jarden \cite[Prop.~6.4.8]{Fried},
states:
\begin{theorem}[Fried--Jarden] \label{friedjarden}
Let $R = \mathbb{F}_q[j]$, let $K$ be its field of fractions, let $L$
be a finite Galois extension of $K$ with Galois group $G$, and let $\mathcal{C} \subset G$ be a conjugacy class.
Let $\mathbb{F}$ be the algebraic
closure of $\mathbb{F}_q$ in $L$ and let $a$ be a positive integer such that $\tau_{\mid_\mathbb{F}}$
acts as $q^a$th power Frobenius for each $\tau \in \mathcal{C}$. Let
$C_1(L/K,\mathcal{C})$ be the set of prime ideals of degree $1$ of $R$ that do
not ramify in $L$ and for which the associated Artin symbol equals $\mathcal{C}$.
If $a \not \equiv 1 \bmod [\mathbb{F} : \mathbb{F}_q]$ then $C_1(L/K, \mathcal{C}) = \emptyset$.
If not, we have
\[ \left| \#C_1(L/K, \mathcal{C}) - \frac{ \# \mathcal{C}}{m} q \right| \, < \, \frac{2 \# \mathcal{C}}{m} \left[ (m + g_L) \sqrt{q} + m \sqrt[4]{q} + g_L + m \right] \]
where $m = [L : K \mathbb{F}]$ and $g_L$ is the genus of $L$ as a function field over $K \mathbb{F}$.
\end{theorem}
For our choice of $R$ and $L$, we have that $\mathbb{F} = \mathbb{F}_q(\zeta_{N'})$.
Then $m = [L : \mathbb{F}_q(\zeta_{N'})(j)]$ is the degree of $J$, i.e.
\[ m = \# \left( \text{SL}_2(Z_{N'} \right) \oplus Z_{p^e}^\times) / \{ \pm 1 \}, \]
and $g_L$ is the genus of $X(p^e; \zeta_{N'})$, which is given by formula (\ref{modulargenus})
above.
An element $\pm (M,u) \in G$ acts as $q$th power Frobenius on $\mathbb{F}_q(\zeta_{N'})$ if and
only if $\det M = q$.

Let $\mathcal{F} \subset \text{GL}_2(Z_{N'})$ and $\mathcal{T} \in Z_{p^e}^\times$
be as in the \'enonc\'e of Theorem~\ref{equidistributionmain}, with the extra condition that we are still assuming $N' > 2$.
Then $\mathcal{F}$ and $\mathcal{T}$ determine a conjugacy class of $\text{GL}_2(Z_{N'}) \oplus Z_{p^e}$
which we abusively denote by $(\mathcal{F}, \mathcal{T})$ and
which is actually contained in $\text{GL}_2^H(Z_{N'}) \oplus Z_{p^e}$. In this smaller group, $(\mathcal{F}, \mathcal{T})$ splits into a union of conjugacy classes
$(\mathcal{F}_1, \mathcal{T}), \dots, (\mathcal{F}_r, \mathcal{T})$ for some $r \in \mathbb{Z}_{\geq 1}$.
Each $(\mathcal{F}_i,\mathcal{T})$ reduces modulo $\{\pm 1\}$ to a conjugacy class $(\overline{\mathcal{F}}_i, \overline{\mathcal{T}})$ of $G$.
Let $(\overline{\mathcal{F}}, \overline{\mathcal{T}})$ denote the union of these conjugacy classes, and let $B = C_1(L/K, (\overline{\mathcal{F}},\overline{\mathcal{T}}))$
be the set of $j_0 \in \mathbb{F}_q$ for which $j_0$ satisfies (\ref{unramifiedcondition}) and
the Artin symbol of $(j-j_0)$ is contained in $(\overline{\mathcal{F}},\overline{\mathcal{T}})$. Then by applying
Theorem~\ref{friedjarden} to each $(\overline{\mathcal{F}}_i,\overline{\mathcal{T}})$ and taking the sum of the resulting
inequalities, we find
\begin{equation} \label{friedjardenappl1}
\left| \#B - \frac{ \# (\overline{\mathcal{F}}, \overline{\mathcal{T}})}{m} q \right| \, < \, \frac{2 \# (\overline{\mathcal{F}}, \overline{\mathcal{T}})}{m} \left[ (m + g_L) \sqrt{q} + m \sqrt[4]{q} + g_L + m \right].
\end{equation}
Now let $A$ denote the set of $\mathbb{F}_q$-isomorphism classes of elliptic curves $E / \mathbb{F}_q$ for
which
\begin{itemize}
\item[-] $j(E) \ \not \in \left\{ \begin{array}{ll} \{0, 1728\} & \text{if $e = 0$} \\
  \mathcal{J}_\text{ss} \cup \{0,1728\} & \text{if $e \geq 1$}; \\ \end{array} \right.$
\item[-] the conjugacy class of $\text{GL}_2(Z_{N'}) \oplus Z_{p^e}^\times$ obtained by
considering the action of $q$th power Frobenius with respect to
all bases $P,Q$ of $E[N']$ equals $\mathcal{F}$;
\item[-] $q$th power Frobenius maps every generator $R$ of $E[p^e]$ to $\mathcal{T} \cdot R$.
\end{itemize}
\begin{lemma} \label{friedjardenappl2}
One has
\[ \left| \#A - \frac{ \# (\mathcal{F}, \mathcal{T})}{m} q \right| \, < \, \frac{4 \# (\mathcal{F}, \mathcal{T})}{m} \left[ (m + g_L) \sqrt{q} + m \sqrt[4]{q} + g_L + m \right]. \]
\end{lemma}
\noindent \textsc{Proof.}
First note that $q$th power Frobenius acts on $E[p^e]$ as multiplication by $\mathcal{T}$ if and only
if it acts on $\text{ker} \, V_{p^e}$ as multiplication by $\mathcal{T}$. By the lemma below, the natural map $A \rightarrow B : E \mapsto j(E)$ is onto and $2$-to-$1$, and for
each $j_0 \in B$, Frobenius acts on the two pre-images with opposite signs.
Thus if $(\mathcal{F}, \mathcal{T}) \cap (- \mathcal{F}, - \mathcal{T}) = \emptyset$, then $\# A = \# B$ and
$\# (\mathcal{F}, \mathcal{T}) = \# (\overline{\mathcal{F}}, \overline{\mathcal{T}})$
and the statement follows from (\ref{friedjardenappl1}).
If $(\mathcal{F}, \mathcal{T}) = (- \mathcal{F}, - \mathcal{T})$ then $\# A = 2 \cdot \# B$ and
$\# (\mathcal{F}, \mathcal{T}) = 2 \cdot \# (\overline{\mathcal{F}}, \overline{\mathcal{T}})$
and the statement again follows. \hfill $\blacksquare$

\begin{lemma} \label{twists}
Let $E/ \mathbb{F}_q$ be an elliptic curve and let $[E]_{\mathbb{F}_q}$
be the set of $\mathbb{F}_q$-isomorphism classes of elliptic curves that are $\overline{\mathbb{F}}_q$-isomorphic
to $E$. Then $\#[E]_{\mathbb{F}_q} \geq 2$. More precisely,
\begin{itemize}
  \item[-] if $j(E) \neq 0, 1728$, then $\#[E]_{\mathbb{F}_q} = 2$ and $[E]_{\mathbb{F}_q}$ consists
  of $E$ and its quadratic twist $E^t$; if $\mathcal{F} \subset \emph{GL}_2(Z_{N'})$ is the conjugacy class determined by $q$th power Frobenius
  acting on $E[N']$, then $-\mathcal{F}$ is the conjugacy class determined by $q$th power Frobenius acting on $E^t[N']$;
  similarly, if $q$th power Frobenius acts on $E[p^e]$ as multiplication
  by $\mathcal{T} \in Z_{p^e}^\times$, then it acts on $E^t[p^e]$ as multiplication by $-\mathcal{T}$;
  \item[-] otherwise, we have the following upper bounds: if $j(E) = 1728$ and $p \neq 2,3$ then $\#[E]_{\mathbb{F}_q} \leq 4$;
  if $j(E) = 0$ and $p \neq 2,3$ then $\#[E]_{\mathbb{F}_q} \leq 6$; if $j(E) = 0 = 1728$ and $p = 3$ then $\#[E]_{\mathbb{F}_q} \leq 12$;
  if $j(E) = 0 = 1728$ and $p = 2$ then $\#[E]_{\mathbb{F}_q} \leq 24$.
\end{itemize}
\end{lemma}
\noindent \textsc{Proof.} We first recall some facts on quadratic twisting, because the existing literature
contains certain ambiguities here. We follow \cite[X.2.4, Exercise~A.2]{Silverman}.
First suppose that $p > 2$. Let
$E$ be an elliptic curve over $\mathbb{F}_q$, take a short Weierstrass model $E: y^2 = f(x)$ and a nonsquare $d \in \mathbb{F}_q$.
Then $E^t$ is defined by $dy^2 = f(x)$. Its $\mathbb{F}_q$-isomorphism class does not depend on the choice of the model, nor on the choice of $d$.
We have an $\overline{\mathbb{F}}_q$-isomorphism $\iota : E^t \rightarrow E : (x,y) \mapsto (x, \sqrt{d}y)$. If $p = 2$ and $j(E) \neq 0$
then $E$ allows a model $y^2 + xy = x^3 + ax^2 + b$ (see \cite[Appendix~A]{Silverman}). Let $d \in \mathbb{F}_q$ have trace $1$, then it is
of the form $\beta^2 + \beta$ for some $\beta \in \mathbb{F}_{q^2} \setminus \mathbb{F}_q$. The quadratic twist $E^t$ is then given by
\[ y^2 + xy = x^3 + (a + d)x^2 + b.\]
This is again well-defined and we have an $\overline{\mathbb{F}}_q$-isomorphism $\iota : E^t \rightarrow E : (x,y) \mapsto (x,y + \beta x)$. Note
that $E$ can a priori be $\mathbb{F}_q$-isomorphic to its quadratic twist, take for instance $q \equiv 3 \bmod 4$, $E: y^2 = x^3 + x$ and $d = - 1$.
Now it is an easy exercise to verify that if $F$ is the matrix of $q$th power Frobenius acting on $E[N']$ with respect to a basis $P,Q$, then
$-F$ is the matrix of $q$th power Frobenius acting on $E^t[N']$ with respect to $\iota^{-1}(P), \iota^{-1}(Q)$, and similarly for $\text{ker} \, V_{p^e}$.

For the remaining statements, we analyze the formula
\[ \sum_{E' \in [E]_{\mathbb{F}_q}} \frac{1}{\# \text{Aut}_{\mathbb{F}_q}(E')} = 1, \]
a proof of which can be found in \cite[Prop.~2.1]{Howe}. Since $\{ \pm 1 \} \subset \text{Aut}_{\mathbb{F}_q}(E')$,
we have that $\# [E]_{\mathbb{F}_q} \geq 2$. The upper bounds follow from $\text{Aut}_{\mathbb{F}_q}(E') \subset \text{Aut}_{\overline{\mathbb{F}}_q}(E')$
and \cite[Thm.~III.10.1]{Silverman}. Finally, if $j(E) \neq 0,1728$, then $E$ cannot be $\mathbb{F}_q$-isomorphic
to its quadratic twist: such an isomorphism would yield a non-rational automorphism of $E$, which cannot exist since $\text{Aut}_{\overline{\mathbb{F}}_q}(E) = \{ \pm 1 \}$. \hfill $\blacksquare$\\

We are now ready to prove Theorem~\ref{equidistributionmain}.

Still assuming $N' > 2$, Lemma~\ref{friedjardenappl2} immediately implies
\[ \left| \# A - \frac{\#\mathcal{F}}{\#\mathcal{M}_{q,N'}} \cdot \frac{1}{\varphi(p^e)} \cdot 2q \right| \leq \frac{8 \cdot \# \mathcal{F}}{\varphi(p^e) \cdot \# \mathcal{M}_{q,N'}} \left[ (m + g_L) \sqrt{q} + m \sqrt[4]{q} + g_L + m \right].\]
By Lemma~\ref{twists}, the number of $\mathbb{F}_q$-isomorphism classes of elliptic curves over $\mathbb{F}_q$ is contained in $[2q, 2q + 22]$.
Taking into account the ramifying $j$-invariants $0$ and $1728$, corresponding to at most $24$ $\mathbb{F}_q$-isomorphism classes, we find
\[ \left| P_{\mathcal{F},\mathcal{T}} - \frac{\#\mathcal{F}}{\#\mathcal{M}_{q,N'}} \cdot \frac{1}{\varphi(p^e)}\right| \qquad \qquad \qquad \qquad \qquad \qquad \qquad \qquad \qquad \qquad \]
\begin{equation} \label{equidistformula}
 \leq \frac{8 \cdot \# \mathcal{F}}{\varphi(p^e) \cdot \# \mathcal{M}_{q,N'}} \cdot \frac{ (m + g_L) \sqrt{q} + m \sqrt[4]{q} + g_L + m + \frac{22}{8}}{2q} + \frac{24}{2q}.
\end{equation}
Note that, in the case $e \geq 1$, both the definition of $A$ and the definition of $P_{\mathcal{F}, \mathcal{T}}$
ruled out all supersingular curves.
Therefore, this has no influence. This is an important difference with the proof of Theorem~\ref{maintheorem} in Section~\ref{traceNdist} below.

Next, we analyze the asymptotical behavior of the error term. From (\ref{modulargenus}) we see that $g_L \leq \frac{p^e \varphi(p^e) N'^3}{24}$, and it is easy
to verify that if $N'$ factors as $\ell_1^{n_1} \cdots \ell_t^{n_t}$ for distinct primes $\ell_i$, then
\[ m = \# \left( \text{SL}_2(Z_{N'}) \oplus Z_{p^e}^\times \right) / \{ \pm 1 \} = \frac{\varphi(p^e)}{2} \cdot \prod_{i = 1}^t \ell_i^{3n_i-2}(\ell_i^2-1) \leq \frac{\varphi(p^e) N'^3}{2}.  \]
Finally,
\[ \# \mathcal{F} \leq \# \{ \text{matrices with trace $\text{Tr}(\mathcal{F})$} \} \leq \prod_{i=1}^t \ell_i^{2n_i - 1}(\ell_i+1) \]
by the results of Section~\ref{traceNdist} below. Hence
\[ \frac{\# \mathcal{F}}{\# \mathcal{M}_{q,N'}} \leq \frac{1}{N'} \prod_{i=1}^t \frac{\ell_i}{\ell_i - 1} \leq \frac{1}{N'} \prod_{\ell \leq 2 \lceil \log_2 N'
\rceil} \frac{\ell}{\ell - 1}\]
where the latter product is over all primes $\ell$. The first inequality follows from the formula for the size of $\mathcal{M}_{q,N'}$ given in Section \ref{traceNdist},  the second inequality follows from
the estimate $\prod_{\ell \leq 2x} \ell \geq 2^x$ (see e.g.\ \cite[Exercise~1.28]{primenumbers}).
Mertens' theorem (see \cite[Corollary~1]{rosser} for an effective version) then shows that
this product is $\mathcal{O}(\ln \ln N')$.
%\[ \frac{\# \mathcal{F}}{\# \mathcal{M}_{q,N'}} \leq \frac{e^\gamma}{N'} \left( 1 + \frac{1}{\ln^2\left( 2 \lceil \log_2 N' \rceil \right)} \right) \ln \left( 2 \lceil \log_2 N' \rceil \right) < \frac{3}{N'} \ln (2\ln N') \]

Theorem~\ref{equidistributionmain} then follows by noting
that the case $N' \leq 2$ is a mere consequence of the case $N' = 4$ (possibly yielding an increase of $C$, though).

\section{The distribution of Frobenius traces} \label{traceNdist}

In this section, we will prove Theorem~\ref{maintheorem} and provide an explicit description
of the function $\Psi$.

For each prime number $p$ and each pair of integers $N'\geq 1$, $e \geq 0$, write $N = N'p^e$ and
define the trace $\text{Tr}(x)$ of an element $x=(M,u)$ of
\[ G = \text{GL}_2(Z_{N'}) \oplus Z_{p^e}^\times \] to be the unique element
of $Z_N$ that is congruent both to $\text{Tr}(M) \bmod N'$ and to $u \bmod p^e$. As before, for each power $q$ of $p$, let $\mathcal{M}_{q,N'} \subset \text{GL}_2(Z_{N'})$ be the set of matrices having determinant $q$. Let $Q$ be as in the introduction and define
\[ r : Q \times \mathbb{Z}_{\geq 1} \times \mathbb{Z}_{\geq 0} \rightarrow \mathbb{Q}_{\geq 0} : (q,N,t) \mapsto
\frac{\# \left\{ \, \left. x \in \mathcal{M}_{q,N'} \times Z_{p^e}^\times \, \right| \, \text{Tr}(x) \equiv t \bmod N \right\} }{\varphi(p^e) \cdot \# \mathcal{M}_{q,N'}}.\]
Then it is easy to verify that $r$ satisfies conditions \emph{(i)} and \emph{(iii)} of Theorem~\ref{maintheorem}.
Moreover, there exists an absolute and explicitly computable constant $C \in \mathbb{R}_{> 0}$ such that, for all $q \in Q, N \in \mathbb{Z}_{\geq 1}$ and $t \in \mathbb{Z}_{\geq 0}$,
\[ \left| P_{q,N}(t) - r(q,N,t) \right| \leq C \cdot \frac{p^eN'^2 \ln \ln N'}{\sqrt{q}}\]
Indeed, let $(\mathcal{F}_t, \mathcal{T}_t)$ denote the set of elements of $\mathcal{M}_{q,N'} \times Z_{p^e}^\times$ having trace $t \bmod N$. It is
a union of conjugacy classes of $G$. By applying Lemma~\ref{friedjardenappl2} to each of these conjugacy classes, taking
the sum of the resulting inequalities, and following a reasoning similar to the one at the end of Section~\ref{matrixdist},
we obtain
\[ \left| P_{q,N}(t) - r(q,N,t) \right| \qquad \qquad \qquad \qquad \qquad \qquad \qquad \qquad \qquad \qquad \qquad \qquad \qquad \]
\begin{equation} \label{tracedistformula}
 \leq \frac{8 \cdot \# (\mathcal{F}_t, \mathcal{T}_t)}{\varphi(p^e) \cdot \# \mathcal{M}_{q,N'}} \cdot \frac{ (m + g_L) \sqrt{q} + m \sqrt[4]{q} + g_L + m + \frac{22}{8}}{2q} + \frac{24}{2q} + \overline{\delta}_{e,0} \frac{p/6}{2q},
\end{equation}
where $\overline{\delta}_{e,0} = 1 - \delta_{e,0}$, with $\delta_{e,0}$ the Kronecker delta. In contrast with $\mathcal{P}_{\mathcal{F}, \mathcal{T}}$, the definition of $P_{q,N}(t)$ does include supersingular curves, whereas the set $A$ from Lemma~\ref{friedjardenappl2} does not as soon as $e \geq 1$.
Therefore, the error term $\frac{24}{2q}$, which in (\ref{equidistformula}) accounted
for the $j$-invariants $0$ and $1728$, should be replaced by an error term accounting in addition
for all supersingular $j$-invariants. For this one can use that there are at most $p/12 + 2$ supersingular
$j$-invariants in $\overline{\mathbb{F}}_q$, and that that for $p=2,3$ the unique supersingular $j$-invariant
is in fact $0=1728$. See \cite[V.4]{Silverman}. An error term analysis as in Section~\ref{matrixdist}
then proves the estimate.

It remains to show that, for the function $\Psi$ that is described below, the function $r$ satisfies condition
\emph{(ii)} of Theorem~\ref{maintheorem}. So assume that $N=\ell^n$ for some integer $n \geq 1$ and a prime
number $\ell$, and let $q \in Q$ be coprime to $\ell$.
Then
\[ r(q,\ell^n,t) = \frac{ \# \left\{ \, M \in \mathcal{M}_{q,\ell^n} \, | \, \text{Tr}(M) \equiv t \bmod \ell^n \, \right\}  }{\# \mathcal{M}_{q,\ell^n}}. \]
It is easy to verify that
$\# \mathcal{M}_{q,\ell^n} = \#\text{SL}_2(Z_{\ell^n})=\ell^{3n-2}(\ell^2-1)$. With $\alpha\in
Z_{\ell^n}\backslash\{0\}$, we define the valuation
$\text{ord}(\alpha)$ as the $\ell$-adic valuation of $\alpha$
embedded in $\ZZ$, whereas we will put $\text{ord}(0)=+\infty$.
Let for $\ell\geq 3$ the map $\Psi:\ZZ\to\ZZ$ be defined as $\Psi = \psi\circ\chi$, where $\chi:\ZZ\to Z_{\ell^n}$ is the natural projection and $\psi: Z_{\ell^n} \rightarrow \mathbb{Z}$ is given by
\[\Delta\mapsto\begin{cases}
\ell^{2n}+\ell^{2n-1} & \text{if $\Delta$ is a nonzero square,}\\
\ell^{2n}+\ell^{2n-1}-2\ell^{2n-\frac k2-1} & \text{if $\Delta$ is no square, $k:=\text{ord}(\Delta)$ is even,}\\
\ell^{2n}+\ell^{2n-1}-(\ell+1)\ell^{2n-\frac{k+3}2} & \text{if $k := \text{ord}(\Delta)$ is odd,}\\
\ell^{2n}+\ell^{2n-1}-\ell^{\frac{3n}2-1} & \text{if $\Delta=0$ and $n$ is even,}\\
\ell^{2n}+\ell^{2n-1}-\ell^{\frac{3n-1}2} & \text{if $\Delta=0$ and $n$ is odd.}\\
\end{cases}
\]
We refer to the end of this section for the definition of
$\Psi$ in case $\ell=2$. Below we prove that indeed:
\begin{theorem}\label{thmTraces}
Let $q$, $t$ and $\ell^n$ be as above and define
$\Delta_t:=t^2-4q$. Then
\[ r(q,\ell^n,t) = \frac{\Psi(\Delta_t)}{\ell^{3n}-\ell^{3n-2}}.\]
\end{theorem}

Let us first discuss some
corollaries. The number of rational points on an elliptic curve
$E$ over $\FF_q$ with trace of Frobenius $T$ equals $q+1-T$. Hence
we can estimate the probability that $\ell^n|\#E(\mathbb{F}_q)$ by
applying Theorem~\ref{thmTraces} with $t=q+1$. Note that then
$t^2-4q\equiv (q-1)^2\bmod \ell^n$. Using this, we recover
the estimates obtained by Howe \cite[Thm.~1.1]{Howe}. %CHH%: new alinea: is a separate corollary

If we suppose $\ell\geq 3$ and $n=1$, then the above formulas
become quite pretty, namely
\[P(t) \sim \begin{cases}\frac{\ell}{\ell^2-1} & \text{if $t^2-4q=0$ in $\FF_\ell$,}\\ \frac{1}{\ell-1} & \text{if $t^2-4q$ is a square in $\FF_l^\times$,}\\ \frac{1}{\ell+1} & \text{if $t^2-4q$ is a nonsquare in $\FF_l$.}\end{cases}\]
The combination of these two corollaries %CHH% the two rather different col's together give Lenstra
generalizes Lenstra's result \cite[Prop.~1.14]{Lenstra} which states that the
probability of $\ell$-torsion approaches $\ell/(\ell^2-1)$ if
$q\equiv 1\bmod \ell$ and $1/(\ell-1)$ otherwise.

The remainder of this section is devoted to the proof of
Theorem \ref{thmTraces}. We note %CHH% We note instead of Note
that the counting of matrices described below was already done by
Gekeler \cite[Thm.~4.4]{Gekeler03} for the case $n\geq
2\cdot \lfloor \frac{\text{ord}(\Delta_t)}{2} \rfloor +2$, using different techniques. %CHH% new alinea

Let $\left(\begin{smallmatrix}u & x\\-y &
z\end{smallmatrix}\right)\in\text{GL}_2(Z_{\ell^n})$ have
determinant $q$ and trace $t$. A trivial computation yields that these conditions are equivalent to the system
of equations
\begin{equation}\label{noncompleted} u=t-z,\qquad xy=z^2-tz+q.\end{equation}
By completing the square, the above system has as many solutions
as
\begin{equation}\label{system}u=t-z,\qquad xy=z^2-\Delta_t/4,\end{equation}
provided that $t/2$ exists modulo $\ell^n$. Suppose for the rest
of the proof that $\ell\geq 3$ and $\Delta_t\in Z_{\ell^n}$; we refer to the end of this
section for the situation $\ell=2$. Clearly all relevant properties (valuation, being a square or not) of $\Delta_t$ and
$\Delta_t/4$ are the same, hence if we can show that the number of
solutions to $xy=z^2-\Delta_t$ equals $\Psi(\Delta_t)$, we are
done. For each value of $z$, we will determine the valuation of
$z^2-\Delta_t$. Then the number of corresponding solutions $(x,y)$
can be computed using the following lemma.
\begin{lemma}\label{lemmaZln}
Let $\ell$ be any prime number, let $n \in \mathbb{Z}_{\geq 1}$
and $\alpha\in Z_{\ell^n}$. Write $k:=\emph{ord}(\alpha)$.
Then the equation $xy=\alpha$ has the following number of
solutions $(x,y)$ in $(Z_{\ell^n})^2$:
\[\begin{cases} (k+1)(\ell^n-\ell^{n-1}) & \text{if $\alpha\not = 0$,}\\ (n+1)(\ell^n-\ell^{n-1})+\ell^{n-1} & \text{if $\alpha=0$}.\end{cases}\]
\end{lemma}
\textsc{Proof.} Suppose $\alpha \neq 0$, the other case works
similarly. We can take $x$ to be any number with valuation
$i\in\{0,1,\ldots,k\}$. For each $i$, the number of such $x$ is
$\ell^{n-i}-\ell^{n-i-1}$. Every choice of $x$ fixes all but the
last $i$ $\ell$-adic digits of $y$, hence we have $\ell^i$
possibilities for $y$. In total this amounts to
\[\sum_{i=0}^k(\ell^{n-i}-\ell^{n-i-1})\ell^i = \sum_{i=0}^k(\ell^n-\ell^{n-1})=(k+1)(\ell^n-\ell^{n-1})\]
solutions $(x,y)$. \hfill$\blacksquare$\\ %CHH% changed squares into blacksquares

Another tool will be the following formula, which is easily proven
by induction:
\begin{lemma} \label{inductionformula}
Let $\ell$ be any prime number, let $n \geq 1$ be an integer and $k \in \{0,1, \dots, n\}$. Then
\[ \sum_{i=0}^k(\ell^{n-i}-\ell^{n-i-1})(2i+1)(\ell^n-\ell^{n-1})= \qquad\qquad\qquad\qquad\qquad\qquad\]
\[\qquad\qquad\qquad\qquad\qquad\qquad\ell^{2n}+\ell^{2n-1}-(2k+3)\ell^{2n-k-1}+(2k+1)\ell^{2n-k-2}.\]
\end{lemma}

Suppose first that $\Delta_t=0$ and $n$ even. Then
$\text{ord}(z^2-\Delta_t)=\text{ord}(z^2)$ for all $z$, and the
number of solutions to $xy = z^2 - \Delta_t$ with
$\text{ord}(z)<n/2$ equals
\[\sum_{i=0}^{n/2-1}(\ell^{n-i}-\ell^{n-i-1})(2i+1)(\ell^n-\ell^{n-1}),\]
by Lemma~\ref{lemmaZln}.
%Clearly the term with index $i$ corresponds to $\text{ord}(z)=i$.
For $\text{ord}(z)\geq n/2$, we find
\[\ell^{n/2}\left((n+1)(\ell^n-\ell^{n-1})+\ell^{n-1}\right)\]
additional solutions. Using Lemma~\ref{inductionformula} one
verifies that the sum of these expressions equals $\Psi(0)$. If
$n$ is odd, then the reasoning is similar.
%, but in this case the first sum stops at $i=\frac{n-1}2$.

Let us now assume that $\Delta_t$ is a nonzero square, i.e.
$\Delta_t=\ell^{2k}\Delta^2$, where $2k<n$ and $\Delta$ is a unit.
Under the change of variables $(x,y,z) \leftarrow (\Delta x,\Delta
y,\Delta z)$ our equation becomes
\begin{equation}\label{equationsimpl}xy=z^2-\ell^{2k}.
%\quad\text{ over $Z_{\ell^n}$.}
\end{equation}
We will use induction on $k$ to show that (\ref{equationsimpl})
has $\Psi(\Delta_t)=\ell^{2n}+\ell^{2n-1}$ solutions. For $k=0$
we have $xy=z^2-1$. If $x$ is any unit, we have $y=x^{-1}(z^2-1)$
and $z$ can be chosen arbitrarily. If $x$ is a nonunit and $y$ is
arbitrary, we have 2 different solutions $z\equiv \pm 1$ modulo
$\ell$, which can both be lifted to $Z_{\ell^n}$. In total this
gives \[(\ell^n-\ell^{n-1})\ell^n+2\ell^{n-1}\ell^n =
\ell^{2n}+\ell^{2n-1}.\] Suppose now that $k\geq 1$.
%In a solution
%of (\ref{equationsimpl}) one of the following two possibilities
%hold:  either $x$ or $y$ is a unit, or both are nonunits.
There are $\ell^{2n}-\ell^{2n-1}$ solutions for which $x$ is a
unit. There are $(\ell^{n}-\ell^{n-1})\ell^{n-1}$ solutions for
which $y$ is a unit and $z$ --- and hence $x$ --- are nonunits. The
solutions for which $x$ and $y$ are both nonunits can be
determined using the induction hypothesis. Indeed, a triplet
$(x,y,z)=(\ell x',\ell y',\ell z')$ satisfies
(\ref{equationsimpl}) if and only if $(x',y',z')$ satisfies
\[x'y'=z'^2-\ell^{2k-2}\quad\text{ over $Z_{\ell^{n-2}}$,}\]
which has $\ell^{2n-4}+\ell^{2n-5}$ solutions.
%in $Z_{\ell^{n-2}}$
For each $x'\in Z_{\ell^{n-2}}$ there are $\ell$ corresponding
values for $x=\ell x'\bmod \ell^n$, and similar for $y$ and $z$.
In total we find then
\[\ell^{2n}-\ell^{2n-1}+(\ell^n-\ell^{n-1})\ell^{n-1}+\ell^3(\ell^{2n-4}+\ell^{2n-5})=\ell^{2n}+\ell^{2n-1}.\]

Next, if $k=\text{ord}(\Delta_t)<+\infty$ is odd, we find the
following sum for the number of solutions
\[\sum_{i=0}^{(k-1)/2}(\ell^{n-i}-\ell^{n-i-1})(2i+1)(\ell-\ell^{n-1}) + \ell^{n-(k+1)/2}(k+1)(\ell^n-\ell^{n-1}),\]
which by Lemma~\ref{inductionformula} equals $\Psi(\Delta_t)$.

Finally, with $k$ even but $\Delta_t$ nonsquare we get
\[\sum_{i=0}^{k/2-1}(\ell^{n-i}-\ell^{n-i-1})(2i+1)(\ell-\ell^{n-1}) + \ell^{n-k/2}(k+1)(\ell^n-\ell^{n-1}),\]
and again the result follows from Lemma~\ref{inductionformula}. This completes the proof for $\ell\geq 3$.\\

We end this section by considering the case $\ell=2$. The appropriate description of $\Psi$ depends now on its argument mod $2^{n+2}$ rather than mod $2^n$. More precisely,
$\Psi = \psi \circ \chi$ where $\chi : \mathbb{Z} \rightarrow Z_{2^{n+2}}$
is the natural projection and $\psi: Z_{2^{n+2}} \rightarrow \mathbb{Z}$ is partially given by
\[\Delta \mapsto \begin{cases}
2^{2n-1} & \text{if $\Delta$ is odd,}\\
2^{2n}+2^{2n-1}-3\cdot 2^{2n-\frac{k+1}2} & \text{if $\Delta\neq 0$ is even and $k:=\text{ord}(\Delta)$ is odd,}\\
2^{2n}+2^{2n-1}-2^{\frac{3n}2-1} & \text{if $\Delta\equiv 0\bmod 2^{n+2}$ and $n$ is even,}\\
2^{2n}+2^{2n-1}-2^{\frac{3n-1}2} & \text{if $\Delta\equiv 0\bmod 2^{n+2}$ and $n$ is odd.}
\end{cases}\]
\noindent In case $\Delta\neq 0$ is even and $\text{ord}(\Delta)=2k>0$ is even as well, the definition of $\psi$ is more complicated. Let $D$
be such that $\Delta = 2^{2k}D$. Then:
\vspace{2mm}

\noindent\begin{tabular}{lrl}
if $n=2k-1$: &  & $\psi(\Delta):=2^{2n}+2^{2n-1} - 2^{\frac{3n-1}2}$,\\
if $n=2k$, & $D\equiv 1\bmod 4$: & $\psi(\Delta):=2^{2n}+2^{2n-1} - 2^{\frac{3n}2-1}$,\\
  & $D\equiv 3\bmod 4$: & $\psi(\Delta):=2^{2n}+2^{2n-1} - 3\cdot 2^{\frac{3n}2-1}$,\\
if $n\geq 2k+1$, & $D\equiv 3\bmod 4$: & $\psi(\Delta):=2^{2n}+2^{2n-1}-3\cdot 2^{2n-k-1}$,\\
 & $D\equiv 5\bmod 8$: & $\psi(\Delta):=2^{2n}+2^{2n-1}- 2^{2n-k}$,\\
& $D\equiv 1\bmod 8$: & $\psi(\Delta):=2^{2n}+2^{2n-1}$.\\
\end{tabular}
\vspace{2mm}

We will now prove that for any $t \in \mathbb{Z}$, the number of
solutions (over $Z_{2^n}$) to the system~(\ref{noncompleted}) is
precisely $\Psi(\Delta_t)$, where $\Delta_t = t^2 - 4q$.
Note first that if $t$ (or equivalently $\Delta_t$) is odd, we have that
$\text{ord}(z^2-tz+q)=0$ for all $z$. Then Lemma~\ref{lemmaZln}
gives a total of $$2^n(2^n-2^{n-1})=2^{2n-1} = \Psi(\Delta_t) $$
solutions.

Therefore suppose that $t$ is even. Then $\Delta_t\equiv 0\bmod
4$, and it makes sense to complete the square in
(\ref{noncompleted}) and analyze the system (\ref{system})
instead. As we are interested in solutions modulo $2^n$, from now
on we will consider $\Delta_t/4$ as an element of $Z_{2^n}$. Note
that this depends on $\Delta_t$ mod $2^{n+2}$. Copying the proofs
of the corresponding cases above, the system (\ref{system}) has
$\Psi(\Delta_t)$ solutions if $\Delta_t / 4 = 0$ (in $Z_{2^n}$)
or if $\text{ord}(\Delta_t/4)<n$ is odd. Hence we assume that
$\text{ord}(\Delta_t / 4) = 2\kappa < n$ is even. Let $D \in
Z_{2^n}$ be such that $2^{2\kappa}D = \Delta_t / 4$. If
$i=\text{ord}(z)<\kappa$ we have $\text{ord}(z^2-\Delta_t/4)=2i$,
so by Lemma~\ref{lemmaZln} and Lemma~\ref{inductionformula} all
such $z$ together account for
\[S := \sum_{i=0}^{\kappa-1}(2^{n-i}-2^{n-i-1})(2i+1)(2^n-2^{n-1}) = 2^{2n}+2^{2n-1}-(2\kappa+3)2^{2n-\kappa-1}\]
solutions $(x,y,z)$. From now on we assume $\text{ord}(z)\geq
\kappa$ and put $z=2^{\kappa}z'$, so that our equation becomes
\[xy=2^{2\kappa}(z'^2-D).\]
Note that $z'$ is only well-determined modulo $2^{n-\kappa}$, and
that we are interested in $z'^2-D\bmod 2^{n-2\kappa}$.

If $n = 2\kappa+1$ we have two possibilities: either $z'\equiv
0\bmod 2$, which gives $2^{n-\kappa-1}(2\kappa+1)2^{n-1}$
solutions $\left(x,y,z' \bmod 2^{n-\kappa}\right)$, or $z'\equiv
1\bmod 2$, which gives $2^{n-\kappa-1}((n+1)2^{n-1}+2^{n-1})$
solutions. If we add $S$ to these two numbers, we find the
requested result.

Let $n = 2\kappa+2$, then we have to distinguish between $D \equiv
1 \bmod 4$ and $D \equiv 3 \bmod 4$. For example, if $D\equiv
3\bmod 4$ and $z'$ is odd, the valuation of $2^{2\kappa}(z'^2-D)$
equals $2\kappa+1$, since $3$ is not a quadratic residue modulo 4.
We leave further details to the reader.

Finally we assume that $n\geq 2\kappa+3$. The cases $D\equiv
3\bmod 4$ and $D\equiv 5\bmod 8$ are similar to the situation
$n=2\kappa+2$ above, so we only go into more details for $D\equiv
1\bmod 8$. Then we know that $D$ is a square modulo
$2^{n-2\kappa}$ and we can proceed as in the case $\ell\geq 3$ and
$\Delta_t$ a nonzero square. However, things work differently for
the induction step $\kappa=0$, i.e.\ $xy=z^2-1\bmod 2^n$, $n\geq
3$. As the valuation of $z^2-1$ cannot be 1 or 2, we have to
consider four situations. Firstly, $\text{ord}(x)=0$, then $z$ can
be chosen arbitrarily and we find $2^{n-1}\cdot 2^n$ solutions.
Secondly, $\text{ord}(x)=1$, then $\text{ord}(y)\geq 2$ and we can
lift the four solutions $z\equiv 1,3,5,7\bmod 8$ to $Z_{2^n}$,
which gives a total of $4\cdot 2^{n-2}2^{n-2}$ solutions. Third,
$\text{ord}(x)=2$ and $\text{ord}(y)\geq 1$ which gives again
$2^{2n-2}$ solutions. Finally, $\text{ord}(x)\geq 3$ and $y$ is
arbitrary, which gives $4\cdot 2^{n-3}2^n$ solutions. Adding all
these terms together gives $2^{2n}+2^{2n-1}$ solutions.

% It is immediate that there are $2^{n-k-1}(2k+1)2^{n-1}$ solutions $(x,y,z')$ with $z'$ even. How many $z'$ give rise to $2^{2k}\text{ord}(z'^2-D)=2k+i$ for $i=1, \ldots, n-2k$? As $z'^2\equiv 1\bmod 8$ as soon as $z'$ is odd, clearly we can take $i\geq 3$. For all these $i<n-2k$ we find $2^{n-k-(i-1)}$ solutions $z'$, hence we have to add
% \[\sum_{i=3}^{n-2k-1}2^{n-k-(i-1)}(2k+i+1)2^{n-1}=(k+2)2^{2n-k-1}+2^{2n-k-2}-(n+2)2^{n+k+1}.\]
% The remaining $2^{k+2}$ values of $z'$ give $2^{2k}(z'^2-D)\equiv 0\bmod 2^n$, so that we indeed find $\varphi(\Delta_t)$ as total amount of solutions.

\section{The probability of a point of order $N$}\label{sectionOrder}

In this section, we prove Theorem~\ref{mainthmorder}. Recall that we defined $P'_q(N)$ as the probability that an elliptic curve over $\mathbb{F}_q$ contains a point of order $N$. %CHH% sentence added

Let $\mathbb{F}_q$ be a finite field of characteristic $p$ with $q$ elements.
Let $E$ be an elliptic curve over $\mathbb{F}_q$. It is well-known (see e.g.\ \cite[Exercise 5.6]{Silverman}) that
\[E(\FF_q) \ \cong\ \ Z_A\oplus Z_B\]
for integers $A,B$ such that $A|B$ and $A|q-1$. Hence if $\text{gcd}(N,q-1)=1$, then $P'_q(N)$ equals the probability $P_{q,N}(q+1)$ that $N|\#E(\FF_q)$,
which can be computed using Theorem~\ref{maintheorem}. However, if $\text{gcd}(N,q-1) > 1$, both probabilities are
fundamentally different.
 The following small example might shed some light on this difference. Let $\ell^n=9$, $q\equiv 1\bmod 9$ and $E$ a random elliptic curve over $\FF_q$. The probability that $\#E(\FF_q)\equiv 0\bmod 9$ approaches (for $q\to\infty$) 11/72. However, the approximate probability that $E$ has a point of order 9 is smaller, namely $9/72$. A corollary is that the probability that $E(\FF_q)[9]\cong Z_3\oplus Z_3$ tends to 2/72.

Entirely analogous to the proof of Theorem~\ref{maintheorem} in Section~\ref{traceNdist}, one sees that
it suffices to consider the case $N = \ell^n$ for some prime $\ell \neq p$ and some integer $n \geq 1$. Moreover, it suffices to prove that the number of matrices in $\text{GL}_2(Z_{\ell^n})$
that are conjugated to a matrix of the form $\left(\begin{smallmatrix} 1 & w\\ 0 & q\end{smallmatrix}\right)$ for a certain $w\in Z_{\ell^n}$
is given by $\theta_{\ell^n} \cdot \# \text{SL}_2(Z_{\ell^n})$, with
\[\theta_{\ell^n}:=\begin{cases} \cfrac{1}{\ell^n-\ell^{n-2}} & \text{if $q\equiv 1\bmod\ell^n$, i.e.\ $\nu\geq n$,} \\  & \\ \cfrac{\ell^{2\nu+1}+1}{\ell^{n+2\nu-1}(\ell^2-1)} & \text{elsewhere},\end{cases}\]
where $\nu$ is the $\ell$-adic valuation of $q-1$. Indeed, $E$ will have an $\mathbb{F}_q$-rational
point of order $\ell^n$ if and only if $\mathcal{F}_E$ is conjugated to an upper diagonal matrix of the above form.

The conjugacy classes of matrices of the form $\left(\begin{smallmatrix}1 & w\\ 0 & q\end{smallmatrix}\right)$ are determined by their representants $M_a$ in Lemma \ref{lemmaMatrixStructure} below. The size of the conjugacy class $\text{Cl}_a$ of $M_a$ can be computed as follows. Let $\text{St}_a$ be the stabilizer subgroup of $M_a$, then the classical orbit-stabilizer theorem states that $\#\text{St}_a\cdot\#\text{Cl}_a=\#\text{GL}_2(Z_{\ell^n})$. Hence it suffices to compute the size of $\text{St}_a$. We know that $\left(\begin{smallmatrix} x & y \\ s & t \end{smallmatrix}\right)\in\text{St}_a$ if and only if $\left(\begin{smallmatrix} x & y \\ s & t \end{smallmatrix}\right)$ is invertible and
\begin{equation}\label{eqStab1}\begin{pmatrix}1 & \ell^a\\ 0 & q\end{pmatrix}\cdot\begin{pmatrix}x & y\\ s & t\end{pmatrix}=\begin{pmatrix}x & y\\ s & t\end{pmatrix}\cdot\begin{pmatrix}1 & \ell^a\\ 0 & q\end{pmatrix}.\end{equation}
This condition is equivalent to the system (using $a\leq \nu$)
\begin{equation}\label{eqStab2}\begin{cases}\ell^as \equiv 0\bmod \ell^n \\ \ell^a(t-x) \equiv y(q-1)\bmod \ell^n.\end{cases}\end{equation}
We can choose $x$ and $y$ at random, so that $t\equiv y(q-1)\ell^{-a}+x\bmod \ell^{n-a}$ and $s\equiv 0\bmod \ell^{n-a}$; we find a total of $\ell^{2n+2a}$ matrices satisfying (\ref{eqStab1}). From these we have to remove the singular matrices, which adds the condition $xt\equiv sy\bmod \ell$. If $a<\nu$ we have by (\ref{eqStab2}) that $s\equiv 0\bmod \ell$ and $t\equiv x\bmod \ell$, hence the only additional restriction is that $x\equiv 0\bmod\ell$. This gives $\ell^{2n+2a-1}$ singular matrices and hence $\#\text{St}_a=\ell^{2n+2a}-\ell^{2n+2a-1}$ for $a<\nu$. If $\nu=n$ it is obvious that $\#\text{Cl}_n=1$, so we are left with considering $\text{St}_\nu$ for $\nu<n$. As shown in the proof of Lemma \ref{lemmaMatrixStructure}, the matrix $\left(\begin{smallmatrix} 1 & \ell^\nu \\ 0 & q \end{smallmatrix}\right)$ is conjugated to $\left(\begin{smallmatrix} 1 & 0 \\ 0 & q \end{smallmatrix}\right)$, and now it is an easy exercise to compute the number $\#\text{St}_\nu = \ell^{2n+2\nu}-(2\ell^{2n-1}-\ell^{2n-2})\ell^{2\nu}$. Combined this gives that the number of matrices conjugated to some $\left(\begin{smallmatrix} 1 & w \\ 0 & q \end{smallmatrix}\right)$ where $\nu<n$ equals (note that $\#\text{GL}_2(Z_{\ell^n})=\ell^{4n-4}(\ell^2-\ell)(\ell^2-1)$):
\[\sum_{a=0}^{\nu-1}\frac{\ell^{4n-4}(\ell^2-\ell)(\ell^2-1)}{\ell^{2n+2a}-\ell^{2n+2a-1}} +
 \frac{\ell^{4n-4}(\ell^2-\ell)(\ell^2-1)}{\ell^{2n+2\nu}-2\ell^{2n+2\nu-1}+\ell^{2n+2\nu-2}}=\ell^{2n}+\ell^{2n-2\nu-1}.
\]
Dividing this number by $\#\text{SL}_2(Z_{\ell^n})$ gives the theorem for $\nu<n$. If $q\equiv 1\bmod\ell^n$ we similarly find
\[\sum_{a=0}^{n-1}\frac{\ell^{4n-4}(\ell^2-\ell)(\ell^2-1)}{\ell^{2n+2a}-\ell^{2n+2a-1}} +
 1=\ell^{2n}.\]
This concludes the proof of Theorem~\ref{mainthmorder}.

\begin{lemma}\label{lemmaMatrixStructure}
Let $\nu=\emph{ord}_\ell(q-1)$. Each matrix over $Z_{\ell^n}$ of the form $\left(\begin{smallmatrix} 1 & w \\ 0 & q \end{smallmatrix}\right)$ is conjugated to precisely one matrix of the set %CHH% added: the set
\[\left\{M_a := \begin{pmatrix} 1 & \ell^a \\ 0 & q\end{pmatrix}\ \right|\left.\vphantom{\begin{pmatrix}1 & \ell^a \\ 0 & q\end{pmatrix}}\ 0\leq a \leq \nu\right\}.\]
\end{lemma}
\textsc{Proof.} First we show that $\left(\begin{smallmatrix} 1 & \ell^a \\ 0 & q \end{smallmatrix}\right)$ with $a\geq \nu$ is conjugated to $\left(\begin{smallmatrix} 1 & \ell^\nu \\ 0 & q \end{smallmatrix}\right)$. Write $q=1+\ell^\nu q'$, then
\[\begin{pmatrix}1 & q'^{-1}(\ell^{a-\nu}-1) \\ 0 & 1 \end{pmatrix}^{-1}\cdot\begin{pmatrix}1 & \ell^a\\ 0 & q\end{pmatrix}\cdot \begin{pmatrix}1 & q'^{-1}(\ell^{a-\nu}-1) \\ 0 & 1 \end{pmatrix}=\begin{pmatrix}1 & \ell^\nu\\ 0 & q\end{pmatrix}.\]
Let $w=\ell^a w'$ with $w'$ a unit in $Z_{\ell^n}$, then
\[\begin{pmatrix}w' & 0 \\ 0 & 1 \end{pmatrix}^{-1}\cdot\begin{pmatrix}1 & \ell^aw'\\ 0 & q\end{pmatrix}\cdot \begin{pmatrix}w' & 0 \\ 0 & 1 \end{pmatrix}=\begin{pmatrix}1 & \ell^a\\ 0 & q\end{pmatrix},\]
which implies that at least one matrix of the above set is conjugated to $\left(\begin{smallmatrix} 1 & w \\ 0 & q \end{smallmatrix}\right)$. The fact that all matrices $M_a$ define different conjugacy classes follows either from a direct reasoning (assuming that two of them are conjugated, the transformation matrix will have determinant 0 modulo $\ell$) or from the computations above which show that the conjugacy classes have different size.\hfill$\blacksquare$\\

\noindent \emph{Note.} It is possible to determine the probability of all kinds of group structures in a similar way. For example, let $0\leq a\leq b$ be integers, $\ell$ a prime coprime to $q$ and suppose we want to know the probability that \[E[\ell^\infty](\FF_q)\cong Z_{\ell^a}\oplus Z_{\ell^b}.\]
This can be done as follows. Let $\mathcal{S}$ be the set of matrices $M$ in $\text{GL}_2(Z_{\ell^{a+b+1}})$ with determinant $q$ for which the following conditions hold:
\begin{itemize}\item[(i)]$\text{Tr}(M)\not\equiv q+1\bmod \ell^{a+b+1}$,
 \item[(ii)]$\text{Tr}(M)\equiv q+1\bmod\ell^{a+b}$,
\item[(iii)] $M$ is conjugated to some $\left(\begin{smallmatrix}1 & w \\ 0 & q\end{smallmatrix}\right)\bmod \ell^b$, and
\item[(iv)] $M\equiv \left(\begin{smallmatrix}1 & 0 \\ 0 & 1\end{smallmatrix}\right)\bmod \ell^{a}$.
\end{itemize}
Then the requested probability tends to $\#\mathcal{S}/\#\text{SL}_2(Z_{\ell^{a+b+1}})$. Note that this question was also
considered by Gekeler in \cite{Gekeler06} in the alternative setting mentioned in the introduction.\\

\noindent \emph{Note.} As pointed out by the anonymous referee of a prior submission of this article,
an alternative proof of Theorem~\ref{mainthmorder} can be obtained by using intermediate results of Howe \cite[Section~4]{Howe}. For each
pair of integers $(M,N)$ for which $M \mid N$, Howe
provides a closed formula for the number of $\mathbb{F}_q$-isomorphism classes of elliptic curves,
counted with a weight that is inversely proportional to the size of the automorphism group, for which
\[ E[N](\mathbb{F}_q) \cong Z_M \oplus Z_N. \]
By letting $M$ range over the divisors of $N$ and summing up the corresponding formulas, one
recovers the estimates from Theorem~\ref{mainthmorder}.

\section*{Acknowledgements}
The authors are very grateful to the anonymous referee of a
prior submission of this document, to Hendrik W.\ Lenstra for suggesting
the use of Chebotarev's density theorem, and to Barry Mazur and Bjorn Poonen
for their helpful comments on modular curves.
Both authors thank F.W.O.-Vlaanderen for its financial support. The first author thanks the Massachusetts
Institute of Technology for its hospitality.

\end{document}